\documentclass[10pt]{amsart}

\usepackage{lipsum}
\usepackage{amsfonts,bm,amssymb,amsmath}
\usepackage{graphicx}
\usepackage{epstopdf}
\usepackage[caption=false]{subfig}
\usepackage{pgfplots}
\usepackage{algorithmic}
\usepackage{amsopn}
\usepackage{amsrefs}
\usepackage{float}
\usepackage{caption}

\usepackage{hyperref}
\usepackage{cleveref}

\usepackage{color}
\usepackage{diagbox}

\usepackage{geometry}

\newtheorem{remark}{Remark}[section]

\hypersetup{
	colorlinks,
	linkcolor={red!50!black},
	citecolor={blue!50!black},
	urlcolor={blue!80!black}
}

\numberwithin{equation}{section}

\definecolor{newcolor1}{rgb}{.8,.349,.1}
\colorlet{bblue}{blue!50!black}
\crefformat{equation}{(#2#1#3)}
\crefmultiformat{equation}{(#2#1#3)}{ and~(#2#1#3)}{, (#2#1#3)}{ and~(#2#1#3)}

\crefformat{figure}{Figure~#2#1#3}
\crefmultiformat{figure}{Figures~ #2#1#3}{ and~#2#1#3}{, (#2#1#3)}{ and~(#2#1#3)}
\crefformat{table}{Table~#2#1#3}

\def\e{\mbox{\boldmath $e$}}
\def\f{\mbox{\boldmath $f$}}

\def\g{\mbox{\boldmath $g$}}
\def\h{\mbox{\boldmath $h$}}
\def\m{\mbox{\boldmath $m$}}

\def\x{\mbox{\boldmath $x$}}
\def\y{\mbox{\boldmath $y$}}

\def\0{\mbox{\boldmath $0$}}

\begin{document}

\title[A third-order method for LLG equation]{Efficient And Stable Third-order Method for Micromagnetics Simulations}

\author[C. Xie]{Changjian Xie}
\address{School of Mathematics and Physics\\ Xi'an-Jiaotong-Liverpool University\\Re'ai Rd. 111, Suzhou, 215123, Jiangsu\\ China.}
\email{Changjian.Xie@xjtlu.edu.cn}

\author[C. Wang]{Cheng Wang}
\address{Mathematics Department\\ University of Massachusetts\\ North Dartmouth\\ MA 02747\\ USA.}
\email{cwang1@umassd.edu}


\date{\today}


\begin{abstract}
To address the magnetization dynamics in ferromagnetic materials described by the Landau-Lifshitz-Gilbert equation under large damping parameters, a third-order accurate numerical scheme is developed by building upon a second-order method  \cite{CaiChenWangXie2022} and leveraging its efficiency. This method boasts two key advantages: first, it only involves solving linear systems with constant coefficients, enabling the use of fast solvers and thus significantly enhancing numerical efficiency over existing first or second-order approaches. Second, it achieves third-order temporal accuracy and fourth-order spatial accuracy, while being unconditionally stable for large damping parameters. Numerical tests in 1D and 3D scenarios confirm both its third-order accuracy and efficiency gains. When large damping parameters are present, the method demonstrates unconditional stability and reproduces physically plausible structures. For domain wall dynamics simulations, it captures the linear relationship between wall velocity and both the damping parameter and external magnetic field, outperforming lower-order methods in this regard.
\end{abstract}

\maketitle

\section{Introduction} 

Given that ferromagnetic materials exhibit bistable intrinsic magnetic order (or magnetization), they are widely used in data storage. To describe the dynamics of this magnetization, researchers rely on the Landau-Lifshitz-Gilbert (LLG) equation \cite{Landau1935On,Gilbert:1955}—a model that incorporates two fundamental dynamic terms: the gyromagnetic term (which conserves energy) and the damping term (which dissipates energy).
Why does the damping term matter? Because it has a direct and strong impact on two critical aspects of magnetic devices: the energy they consume and the speed at which they operate. Notably, a recent experiment on magnetic-semiconductor heterostructures \cite{Zhang2020ExtremelyLM} demonstrated that the Gilbert damping constant—a key parameter in the damping term—can be adjusted. At the microscopic scale, damping is driven by three main mechanisms: electron scattering, itinerant electron relaxation \cite{Heinrich1967TheIO}, and phonon-magnon coupling \cite{Suhl1998TheoryOT, Nan2020ElectricfieldCO}. Importantly, these mechanisms can be quantified through electronic structure calculations \cite{TangXia2017}, which is valuable for practical applications. For instance, by tuning the damping parameter, engineers can optimize a material’s magnetodynamic properties—such as reducing the switching current required for magnetic memory devices or increasing their writing speed \cite{Wei2012MicromagneticsAR}.
While most experimental studies have concentrated on scenarios with small damping parameters \cite{Budhathoki2020LowGD,Lattery2018LowGD,Weber2019GilbertDO}, there is growing evidence of large damping effects in specific cases. For example, \cite{GilbertKelly1955, Tanaka2014MicrowaveAssistedMR} document such effects, with \cite{Tanaka2014MicrowaveAssistedMR} further noting that a larger damping constant leads to a shorter magnetization switching time. In fact, \cite{GilbertKelly1955} reports extremely large damping parameters (on the order of 9), highlighting the need to address this regime.

The LLG equation itself poses unique challenges: it is a vectorial, nonlinear system, and it enforces that magnetization length remains constant at every point. To tackle these challenges, significant research has focused on developing numerical methods for micromagnetics simulations that are both efficient and stable (for comprehensive reviews, see \cite{kruzik2006recent,cimrak2007survey}). Among the most widely used approaches are semi-implicit schemes, which offer a key advantage: they avoid the need for complex nonlinear solvers while still maintaining numerical stability \cite{alouges2006convergence, gao2014optimal, Xie2018}.
To illustrate, \cite{Xie2018} developed a second-order backward differentiation formula (BDF) scheme using one-sided interpolation. However, this scheme has a limitation: at each time step, it requires solving a three-dimensional linear system with non-constant coefficients. Later, \cite{jingrun2019analysis} addressed this by establishing a theoretical framework that confirms the scheme’s second-order convergence. Another semi-implicit approach, proposed in \cite{alouges2006convergence}, uses the tangent space to enforce the magnetization length constraint—though it only achieves first-order temporal accuracy. More recently, \cite{Lubich2021} expanded on this work by constructing and analyzing high-order BDF schemes.
Despite these advances, semi-implicit schemes have a critical drawback. While \cite{jingrun2019analysis,Lubich2021} proved that these schemes have unconditional unique solvability, their convergence analysis requires that the temporal step-size is proportional to the spatial grid-size. Even more problematic: the LLG equation’s vectorial structure results in non-symmetric linear systems at each time step. This means FFT-based fast solvers— which are highly efficient for symmetric systems—cannot be used. Instead, researchers often rely on GMRES, but its efficiency is heavily dependent on both the temporal step-size and spatial grid-size. Worse, extensive numerical experiments have shown that GMRES is significantly more computationally costly than standard Poisson solvers \cite{Xie2018}.
To address these limitations, this paper proposes a new high-order numerical method for solving the LLG equation in the presence of large damping parameters. The method achieves third-order accuracy in time and fourth-order accuracy in space, while its computational complexity is comparable to that of solving the scalar heat equation. How is this accomplished? First, the LLG system is reformulated: the damping term is re-expressed as a harmonic mapping flow. Next, the constant-coefficient Laplacian component is discretized using a standard BDF3 temporal scheme—its associated dissipation provides the foundation for numerical stability. Finally, all nonlinear components (including both the gyromagnetic term and the remaining nonlinear expansions in the damping term) are approximated using a fully explicit third-order extrapolation formula.
This explicit treatment of nonlinear terms is a game-changer: it means the scheme only requires a standard Poisson solver at each time step. Because the linear system involved has a symmetric positive definite (SPD) structure, FFT-based fast solvers can be efficiently applied—dramatically reducing computational effort. Extensive numerical experiments not only confirm the scheme’s stability but also reveal important insights: pre-projected solutions introduce significant instability, even though the dissipative property of the heat equation component can partially ensure the stability of nonlinear parts under large damping parameters. Additionally, higher-order numerical methods (e.g., BDF2, BDF3) are less stable than lower-order ones (e.g., BDF1). This is because higher-order methods depend more heavily on solution information from previous unprojected steps, which has a greater negative impact on stability—particularly when updating the stray field.

The rest of the paper is organized to build on this foundation: \cref{sec: numerical scheme} first reviews the micromagnetics model, then presents the proposed numerical method, and finally compares it to the BDF1 and BDF2 methods. \cref{sec:experiments} follows with numerical results, including checks of temporal and spatial accuracy (in both 1D and 3D computations), investigations of numerical efficiency (compared to BDF1 and BDF2), studies of stability with respect to the damping parameter, and an analysis of domain wall motion instability. Finally, \cref{sec:conclusions} offers concluding remarks.

\section{The physical model and the numerical method}
\label{sec: numerical scheme}

\subsection{Landau-Lifshitz-Gilbert equation}

The LLG equation describes the dynamics of magnetization which consists of the gyromagnetic term and the damping term~\cite{Landau1935On,Brown1963micromagnetics}. In the nondimensionalized form, this equation reads as
\begin{align}\label{eq-5}
\m_t=-\m\times(\epsilon\Delta\m+\f)-\alpha\m\times\m\times(\epsilon\Delta\m+\f),
\end{align}
where 
the following source term is defined 
\begin{align}\label{eq-4}
\f=-q(m_2\e_2+m_3\e_3)+\h_s+\h_e,
\end{align}
where $\h_e$ is the external field and $\h_s$ is the stray field with the following formula 
\begin{align}\label{eqn:div}
{\h}_{\text{s}}=\frac{1}{4\pi}\nabla \int_{\Omega} \nabla\left( \frac{1}{|\x-\y|}\right)\cdot {\bm m}({\bm y})\,d{\bm y}.
\end{align}
Here, the dimensionless parameters become $\epsilon=C_{ex}/(\mu_0 M_s^2L^2)$ and $q=K_u/(\mu_0 M_s^2)$ with $L$ the diameter of the ferromagnetic body and $\mu_0$ the permeability of vacuum. The unit vectors are given by ${\bm e}_2=(0,1,0)$, ${\bm e}_3=(0,0,1)$, and $\Delta$ denotes the standard Laplacian operator. For the Permalloy, an alloy of Nickel ($80\%$) and Iron ($20\%$), typical values of the physical parameters are given by: the exchange constant $C_{ex}=1.3\times 10^{-11}\,\textrm{J/m}$, the anisotropy constant $K_u = 100\, \textrm{J/}\textrm{m}^3$, the saturation magnetization constant $M_s = 8.0\times 10^{5}\,\textrm{A/m}$. If $\Omega$ is a rectangular domain, the evaluation of \eqref{eqn:div} can be efficiently done by the Fast Fourier Transform (FFT) \cite{Wang2000}.

Thanks to point-wise identity $|\m|=1$, we obtain an equivalent form: 
\begin{equation}\label{eq-model}
\m_t=\alpha  (\epsilon\Delta\m+\f)+\alpha \left(\epsilon |\nabla \m|^2 -\m \cdot\f \right)\m-\m\times(\epsilon\Delta\m+\f).
\end{equation}
with the homogeneous Neumann boundary condition
\begin{equation}\label{boundary-large}
\frac{\partial{\m}}{\partial {\bm \nu}}\Big|_{\partial \Omega}=0,
\end{equation}
where $\Omega$ is a bounded domain occupied by the ferromagnetic material and $\bm \nu$ is unit outward normal vector along $\partial \Omega$. 

In more details, the magnetization ${\m}\,:\,\Omega\subset\mathbb{R}^d\to \mathbb{R}^3,d=1,2,3 $ is a three-dimensional vector field with a pointwise constraint $|\m|=1$. The first term on the right-hand side in \cref{eq-5} is the gyromagnetic term and the second term stands for the damping term, with $\alpha>0$ being the dimensionless damping coefficient.

In particular, it is noticed that the damping term is rewritten as a harmonic mapping flow, which contains a constant-coefficient Laplacian diffusion term. This fact will greatly improve the numerical stability of the proposed scheme.  

For the numerical description, we first introduce some notations for discretization and numerical approximation. 
Denote the temporal step-size by $k$, and $t^n=nk$, $n\leq \left\lfloor\frac{T}{k}\right\rfloor$ with $T$ the final time. The spatial mesh-size is given by $h_x=h_y=h_z=h=1/N$, and $\m_{i,j,\ell}^n$ stands for the magnetization at time step $t^n$, evaluated at the spatial location $(x_{i-\frac12},y_{j-\frac12},z_{\ell-\frac12})$ with $x_{i-\frac12}=\left(i-\frac12\right)h_x$, $y_{j-\frac12}=\left(j-\frac12\right)h_y$ and $z_{\ell-\frac12}=\left(\ell-\frac12\right)h_z$ ($0\leq i,j,\ell\leq N+1$). In addition, a third order extrapolation formula is used to approximate the homogeneous Neumann boundary condition. For example, such a formula near the boundary along the $z$ direction is given by 
\begin{align*}
\m_{i,j,1}=\m_{i,j,0},\quad \m_{i,j,-1}=\m_{i,j,2},\quad \m_{i,j,N+1}=\m_{i,j,N}\quad \m_{i,j,N+2}=\m_{i,j,N-1}.
\end{align*}
The boundary extrapolation along other boundary sections can be similarly made.

%

The standard second-order centered difference applied to $\Delta \m$ results in
\[
\Delta_h \m_{i,j,k} = \frac{\delta_x^2 \m_{i,j,k}}{h_x^2} + \frac{\delta_y^2 \m_{i,j,k}}{h_y^2} + \frac{\delta_z^2 \m_{i,j,k}}{h_z^2},
\]
where $\delta_x^2 \m_{i,j,k} = \m_{i+1,j,k} - 2\m_{i,j,k} + \m_{i-1,j,k}$, $\delta_y^2 \m_{i,j,k} = \m_{i,j+1,k} - 2\m_{i,j,k} + \m_{i,j-1,k}$ and $\delta_z^2 \m_{i,j,k} = \m_{i,j,k+1} - 2\m_{i,j,k} + \m_{i,j,k-1}$ and the discrete gradient operator $\nabla_h \m$ with $\m = (u, v, w)^T$ reads as
\[
\nabla_h \m_{i,j,k} = \left( \frac{\delta_x \m_{i,j,k}}{h_x}, \frac{\delta_y \m_{i,j,k}}{h_y}, \frac{\delta_z \m_{i,j,k}}{h_z} \right)^T,
\]
where $\delta_x \m_{i,j,k} = \m_{i+1,j,k} - \m_{i-1,j,k}$, $\delta_y \m_{i,j,k} = \m_{i,j+1,k} - \m_{i,j-1,k}$ and $\delta_z \m_{i,j,k} = \m_{i,j,k+1} - \m_{i,j,k-1}$.

Subsequently, the BDF1 and the BDF2 numerical methods need to be reviewed, which could be used for the later comparison.

\subsection{The first order method}

The first-order BDF (BDF1) method is based on a Backward
Differentiation Formula, combined with an explicit extrapolation. The cross product nonlinear
term is treated with previous step solution. It only requires a linear equation solvers with constant
coefficients; as a result, the FFT-based fast solvers could be easily applied. This method is first-order in time and second-order in space. Below is the detailed outline of the BDF1 method.

%

\begin{equation}\label{BDF1}
\left\{ 
\begin{aligned}
&\frac{\tilde{\m}_h^{n+1} - {\m}_h^{n} }{k}
=  - {\m}_h^{n} \times \left(\epsilon \Delta_h {{\m}}_h^{n} +{\f}_h^{n}\right)  + \alpha \left(\epsilon \Delta_h \tilde{\m}_h^{n+1}+{\f}_h^{n}\right)\\  
&\quad + \alpha \left(\epsilon | \nabla_h {\m}_h^{n} |^2-{\m}_h^{n}\cdot {\f}_h^{n}\right) {\m}_h^{n},\\ 
&\m_h^{n+1} = \frac{\tilde{\m}_h^{n+1}}{ |\tilde{\m}_h^{n+1}| } ,
\end{aligned}
\right.
\end{equation}


\subsection{The second-order method}
Such approach has been outlined in \cite{CaiChenWangXie2022}. This method is
based on the second-order BDF temporal discretization, combined with an explicit extrapolation.
It is found that BDF2 is unconditionally stable and is second-order accurate in both space and time. The algorithmic details are given as follows.
\begin{equation}\label{sipm}
\left\{ 
\begin{aligned}
&\frac{\frac32 \tilde{\m}_h^{n+2} - 2 {\m}_h^{n+1} + \frac12 {\m}_h^n}{k}
=  - \hat{\m}_h^{n+2} \times \left(\epsilon \Delta_h \hat{{\m}}_h^{n+2} +\hat{\f}_h^{n+2}\right) \\
&\quad + \alpha \left(\epsilon \Delta_h \tilde{\m}_h^{n+2}+\hat{\f}_h^{n+2}\right)\\  
&\quad + \alpha \left(\epsilon | \nabla_h \hat{\m}_h^{n+2} |^2-\hat{\m}_h^{n+2}\cdot \hat{\f}_h^{n+2}\right) \hat{\m}_h^{n+2},\\ 
&\m_h^{n+2} = \frac{\tilde{\m}_h^{n+2}}{ |\tilde{\m}_h^{n+2}| } ,
\end{aligned}
\right.
\end{equation}
where $\tilde{\m}_h^{n+2}$ is an intermediate magnetization, and $\hat{\m}_h^{n+2}$, $\hat{\f}_h^{n+2}$ are given by the following extrapolation formula: 
\begin{align*}
\hat{\m}_h^{n+2} &=2{\m}_h^{n+1}-{\m}_h^n, \label{m_hat}\\
\hat{\f}_h^{n+2} &=2{\f}_h^{n+1}-{\f}_h^n,
\end{align*}
with $\f_h^{n}=-Q(m_2^n\e_2+m_3^n\e_3)+\h_s^n+\h_e^n$. 

\subsection{The proposed third order method} \label{discretisations}

The high-order BDF idea leads to the proposed numerical method as follows.
\begin{equation}\label{proposed}
\left\{ 
\begin{aligned}
&\frac{\frac{11}{6} \tilde{\m}_h^{n+3} - 3 {\m}_h^{n+2} + \frac32 {\m}_h^{n+1}-\frac13 {\m}_h^n}{k}
=  - \hat{\m}_h^{n+3} \times \left(\epsilon \Delta_{h,(4)} \hat{{\m}}_h^{n+3} +\hat{\f}_h^{n+3}\right) \\
&\quad + \alpha \left(\epsilon \Delta_h \tilde{\m}_h^{n+3}+\hat{\f}_h^{n+3}\right)\\  
&\quad + \alpha \left(\epsilon | \nabla_h \hat{\m}_h^{n+3} |^2-\hat{\m}_h^{n+3}\cdot \hat{\f}_h^{n+3}\right) \hat{\m}_h^{n+3},\\ 
&\m_h^{n+3} = \frac{\tilde{\m}_h^{n+3}}{ |\tilde{\m}_h^{n+3}| } ,
\end{aligned}
\right.
\end{equation}
where
\begin{align*}
\hat{\m}_h^{n+3} &= 3 \m_h^{n+2}-3\m_h^{n+1} + \m_h^n,\\
\hat{\f}_h^{n+3} &= 3 \f_h^{n+2} - 3\f_h^{n+1}+\f_h^n.
\end{align*}

\cref{tab-features} compares the proposed method, the BDF2 and the BDF1 in terms of number of unknowns, dimensional size, symmetry pattern, and availability of FFT-based fast solver of linear systems of equations, and the number of stray field updates. At the formal level, the proposed method is clearly superior to both the BDF2 and the BDF1 algorithms. 
In
 more details, this scheme will greatly improve the computational efficiency, since only three Poisson
solvers are needed at each time step. Moreover, this numerical method preserves a third-order accuracy in time and fourth-order accuracy in space. Interestingly, the numerical results in \cref{sec:experiments} will
demonstrate that the proposed scheme provides a subtle approach for micromagnetics simulations
with less stability in the regime of large damping parameters.

\begin{table}[htbp]
	\begin{center}
		\caption{Comparison of the proposed method, the BDF2 method, and the BDF1 method.}\label{tab-features}
		\begin{tabular}{cccc}
			\hline
			Property or number & Proposed method & BDF2 & BDF1\\
			\hline
			Linear systems& \boldsymbol{$3$} & $3$ & $3$ \\
			Size & \boldsymbol{$N^3$} & $N^3$& $N^3$ \\
			Symmetry& {\bf Yes}& Yes& Yes \\
			Fast Solver& {\bf Yes}& Yes& Yes \\
			Accuracy& \boldsymbol{$\mathcal{O}(k^3+h^4)$} & $\mathcal{O}(k^2+h^2)$ & $\mathcal{O}(k+h^2)$ \\
			Stray field updates & \boldsymbol{$1$} &$1$ &$1$ \\
			\hline
		\end{tabular}
	\end{center}
\end{table}

\begin{remark}
To kick start the proposed method, one can apply a first-order and a second-order
algorithm, such as the first-order BDF method and the second-order BDF method, in the first and
second time step. An overall third-order accuracy is preserved in such an approach.
\end{remark}

\section{Numerical experiments}
\label{sec:experiments}

In this section, we present a few numerical experiments with a sequence of damping parameters for the proposed method, the BDF2 \cite{CaiChenWangXie2022} and the BDF1 method, with the accuracy, efficiency, and stability examined in details. Domain wall dynamics is studied and its velocity is recorded in terms of the damping parameter and the external magnetic field. 

\subsection{Accuracy and efficiency tests}

We set $\epsilon=1$ and $\f=0$ in \cref{eq-model} for convenience. The 1D exact solution is given by 
\begin{equation*}
	\m_e=\left(\cos(\cos(\pi x))\sin t, \sin(\cos(\pi x))\sin t, \cos t\right)^T,
\end{equation*}
and the corresponding exact solution in 3D becomes 
\begin{equation*}
\m_e=\left(\cos(\cos(\pi x)\cos(\pi y)\cos(\pi z))\sin t, \sin(\cos(\pi x)\cos(\pi y)\cos(\pi z))\sin t, \cos t\right)^T.
\end{equation*}
In fact, the above exact solutions satisfy \cref{eq-model} with the forcing term $\g=\partial_t \m_e-\alpha \Delta \m_e -\alpha |\nabla \m_e|^2+\m_e \times \Delta \m_e$, as well as the homogeneous Neumann boundary condition. 

For the temporal accuracy test in the 1D case, we fix the spatial resolution as $h=1D-4$, so that the spatial approximation error becomes negligible. The damping parameter is taken as $\alpha=10$, and the final time is set as $T=0.1$. In the 3D test for the temporal accuracy, 
due to the limitation of spatial resolution, we take a sequence of spatial and temporal mesh sizes: $k=h_x^2=h_y^2=h_z^2=h^2=T/N_0$ for the first-order method and $k=h_x=h_y=h_z=h=T/N_0$ for the second-order method, and $k=h_x^{\frac43}=h_y^{\frac43}=h_z^{\frac43}=h^{\frac43}=T/N_0$ for the proposed method, with the variation of $N_0$ indicated below. Similarly, the damping parameter is given by $\alpha=10$, while the final time $T$ is indicated below. In turn, the numerical errors are recorded in term of the temporal step-size $k$ in \cref{tab-1}. It is clear that the temporal accuracy orders of the proposed numerical method, the BDF2, and the BDF1 are given by $3$, $2$, and $1$, respectively, in both the 1D and 3D computations. 

\begin{table}[htbp] 
	\centering
	{\caption{The numerical errors for the proposed method, the BDF1 and the BDF2 with $\alpha=10$ and $T=0.1$. Left: 1D with $h=1D-4$; Right: 3D with $k=h_x^2=h_y^2=h_z^2=h^2=T/N_0$ for BDF1 and $k=h_x=h_y=h_z=h=T/N_0$ for the BDF2, and $k=h_x^{\frac43}=h_y^{\frac43}=h_z^{\frac43}=h^{\frac43}=T/N_0$ for the proposed method, with $N_0$ specified in the table.}\label{tab-1} }{
		\subfloat[Proposed method]{\label{tab:floatrow:one}%
			\begin{tabular}{cccc|cccc} 
				\hline	
				1D  & & {} & {} &3D &{} & {} &{} \\
				$k$ & $\|\cdot\|_{\infty}$ & $\|\cdot\|_{2}$ & $\|\cdot\|_{H^1}$ & 	$k,k^3\approx h^4$ & $\|\cdot\|_{\infty}$ & $\|\cdot\|_{2}$ & $\|\cdot\|_{H^1}$ \\
				\hline
			$T/8$& 1.977D-7& 1.286D-7& 5.585D-7 & $T/4$ & 8.548D-6 & 1.758D-6 & 3.557D-5 \\	
			$T/12$ &5.858D-8 & 3.826D-8 & 1.660D-7 & $T/5$ & 3.924D-6 & 8.345D-7 & 1.833D-5 \\
			$T/16$ & 2.461D-8 & 1.612D-8 & 6.994D-8&$T/6$ & 2.257D-6 & 4.805D-7 &  1.008D-5\\
			$T/24$ & 7.254D-9 & 4.767D-9 & 2.070D-8 &$T/8$ & 1.001D-6 & 2.103D-7 & 4.481D-6\\
			$T/32$ & 3.055D-9 & 2.009D-9 & 8.728D-9& $T/9$& 7.546D-7 & 1.584D-7 & 3.362D-6\\
				order &3.01 &3.00& 3.00&{--}&2.98&2.97&2.93\\
				
				
				%
				\hline
			\end{tabular}	
		}
		\qquad
		\subfloat[BDF1]{\label{tab:floatrow:two}
			\begin{tabular}{cccc|cccc} 
				\hline
				1D & &  & {} & 3D & {} &  & {} \\
				$k$ & $\|\cdot\|_{\infty}$ & $\|\cdot\|_{2}$ & $\|\cdot\|_{H^1}$ & $k=h^2$ & $\|\cdot\|_{\infty}$ & $\|\cdot\|_{2}$ & $\|\cdot\|_{H^1}$\\
				\hline
				$T/8$ & 1.267D-3 & 8.289D-4 & 3.613D-3& $T/40$ & 5.381D-4 & 8.166D-5 & 5.912D-4 \\
				$T/12$ &8.420D-4 & 5.507D-4 & 2.396D-3 & $T/57$ & 3.764D-4 & 5.682D-5 & 4.109D-4 \\
				$T/16$ &6.301D-4 & 4.124D-4 & 1.793D-3 & $T/78$ & 2.763D-4 & 4.158D-5 & 3.009D-4  \\
				$T/24$ & 4.188D-4 & 2.745D-4 & 1.193D-3 &$T/102$&2.116D-4 & 3.178D-5 & 2.300D-4 \\
				$T/32$ & 3.133D-4 & 2.057D-4 & 8.940D-4 &$T/129$ & 1.674D-4 & 2.510D-5 & 1.815D-4\\
				order &1.01 &1.01&1.01 &{--}&1.00 &1.01&1.01 \\
				
				
				%
				\hline
			\end{tabular}
		}
		\qquad
		\subfloat[BDF2]{\label{tab:floatrow:three}
			\begin{tabular}{cccc|cccc} 
				\hline
				1D  & & & {} & 3D & {} &  & {}\\
				$k$ & $\|\cdot\|_{\infty}$ & $\|\cdot\|_{2}$ & $\|\cdot\|_{H^1}$ & $k=h$ & $\|\cdot\|_{\infty}$ & $\|\cdot\|_{2}$ & $\|\cdot\|_{H^1}$\\
				\hline
				$T/8$ &5.240D-6& 3.700D-6& 9.774D-6& $T/3$ &1.639D-4 & 2.824D-5 & 2.272D-4\\	
			$T/12$ &2.468D-6& 1.748D-6& 4.501D-6& $T/4$ &9.426D-5 & 1.619D-5 & 1.293D-4 \\
				$T/16$ & 1.428D-6 & 1.013D-6 & 2.576D-6 &$T/5$ &6.150D-5 & 1.060D-5 & 8.266D-5 \\
			$T/24$ &6.527D-7& 4.636D-7& 1.165D-6&	 $T/6$ &4.320D-5 & 7.490D-6 & 5.741D-5 \\
				$T/32$ & 3.725D-7& 2.646D-7& 6.615D-7 & $T/7$& 3.201D-5 & 5.570D-6 & 4.221D-5 \\
				order &1.91&1.90&1.94 & {--}&1.93 &1.91&1.99\\
				
				
				%
				\hline
			\end{tabular}
	} }	
\end{table}

The spatial accuracy order is tested by fixing $k=1D-5$, $\alpha=10$, $T=0.1$ in 1D and $k=1D-4$, $\alpha=10$, $T=0.1$ in 3D. The numerical error is recorded in term of the spatial grid-size $h$ in \cref{tab-2}. Similarly, the presented results have indicated the fourth order spatial accuracy of all the proposed algorithms and the second order spatial accuracy of the BDF2, and the BDF1 methods, in both the 1D and 3D computations. 

\begin{table}[htbp]
	\centering
	{\caption{The numerical errors of the proposed method, the BDF1 and the BDF2 with $\alpha=10$ and $T=0.1$. Left: 1D with $k=1D-5$; Right: 3D with $k=1D-4$.} \label{tab-2} }{
		\subfloat[Proposed method]{\label{tab:floatrow:1-S}
			\begin{tabular}{cccc|cccc}	
				\hline
				1D  & & & {} & 3D & & & \\
				$h$ & $\|\cdot\|_{\infty}$ &$\|\cdot\|_{2}$ &$\|\cdot\|_{H^1}$& $h$ & $\|\cdot\|_{\infty}$ & $\|\cdot\|_{2}$ & $\|\cdot\|_{H^1}$  \\
				\hline
				$1/16$ &7.525D-6 &5.707D-6 &9.903D-5 & $1/16$ & 7.064D-6& 1.689D-6& 3.566D-5\\
				$1/32$ &4.916D-7& 3.625D-7& 6.383D-6 & $1/20$ & 2.989D-6 & 7.009D-7 & 1.469D-5\\
				$1/64$ & 3.108D-8& 2.276D-8& 4.021D-7 & $1/24$ & 1.465D-6 & 3.434D-7 & 7.096D-6\\
				$1/128$ & 1.948D-9& 1.424D-9& 2.519D-8 & $1/28$ & 7.975D-7 & 1.897D-7 & 3.827D-6 \\
				$1/256$ & 1.215D-10& 8.899D-11 & 1.575D-9 & $1/32$ & 4.689D-7 & 1.153D-7 & 2.239D-6 \\
				order  &3.98 &3.99& 3.99& {--} & 3.91&3.88&3.99  \\
				\hline
			\end{tabular}
		}	
		\qquad
		\subfloat[BDF1]{\label{tab:floatrow:2-S}
			\begin{tabular}{cccc|cccc}	
				\hline
				1D  & &  & {} & 3D & & & \\
				$h$ & $\|\cdot\|_{\infty}$ &$\|\cdot\|_{2}$ &$\|\cdot\|_{H^1}$ & $h$ & $\|\cdot\|_{\infty}$ &$\|\cdot\|_{2}$ &$\|\cdot\|_{H^1}$\\
				\hline
				$1/16$ & 3.069D-4 & 2.861D-4 & 2.000D-3 & $1/16$ &4.562D-4 & 9.178D-5 & 7.985D-4\\
				$1/32$ & 7.756D-5 & 7.122D-5 & 4.988D-4 & $1/20$ & 2.942D-4 & 5.864D-5 & 5.076D-4 \\
				$1/64$ & 1.943D-5 & 1.779D-5 & 1.246D-4 & $1/24$ & 2.053D-4 & 4.069D-5 & 3.509D-4 \\
				$1/128$ & 4.851D-6 & 4.449D-6 & 3.115D-5 & $1/28$ &1.514D-4 & 2.988D-5 & 2.570D-4 \\
				$1/256$ & 1.204D-6 & 1.117D-6 & 7.788D-6 & $1/32$ &1.163D-4 & 2.288D-5 & 1.963D-4 \\
				order  &2.00 &2.00&2.00 & {--} & 1.97&2.01&2.02 \\
				\hline
			\end{tabular}
		}
		\qquad
		\subfloat[BDF2]{\label{tab:floatrow:3-S}
			\begin{tabular}{cccc|cccc}	
				\hline
				1D  & &  & & 3D & & & \\
				$h$ & $\|\cdot\|_{\infty}$ &$\|\cdot\|_{2}$ &$\|\cdot\|_{H^1}$ &$h$ & $\|\cdot\|_{\infty}$ &$\|\cdot\|_{2}$ &$\|\cdot\|_{H^1}$ \\
				\hline
				$1/16$ &3.069D-4 & 2.861D-4 & 2.000D-3 & $1/16$ & 4.553D-4& 9.175D-5& 7.986D-4\\
				$1/32$ &7.757D-5 & 7.122D-5 & 4.988D-4 & $1/20$ & 2.934D-4 & 5.860D-5& 5.077D-4 \\
				$1/64$ & 1.944D-5 & 1.778D-5 & 1.246D-4 & $1/24$ & 2.045D-4 & 4.065D-5 & 3.510D-4\\
				$1/128$ &4.863D-6 & 4.445D-6 & 3.115D-5 & $1/28$ &1.506D-4 & 2.985D-5 & 2.571D-4 \\
				$1/256$ &1.216D-6 & 1.111D-6 & 7.788D-6 & $1/32$ & 1.154D-4 & 2.284D-5 & 1.964D-4 \\
				order  &2.00 &2.00&2.00& {--} &1.98&2.01&2.02 \\
				\hline
			\end{tabular} 
		}
	}
\end{table}

\begin{figure}[htbp]
	\centering
	\subfloat[Varying $k$ in 1D up to $T=1$ ]{\label{cputime_1D}\includegraphics[width=2.5in]{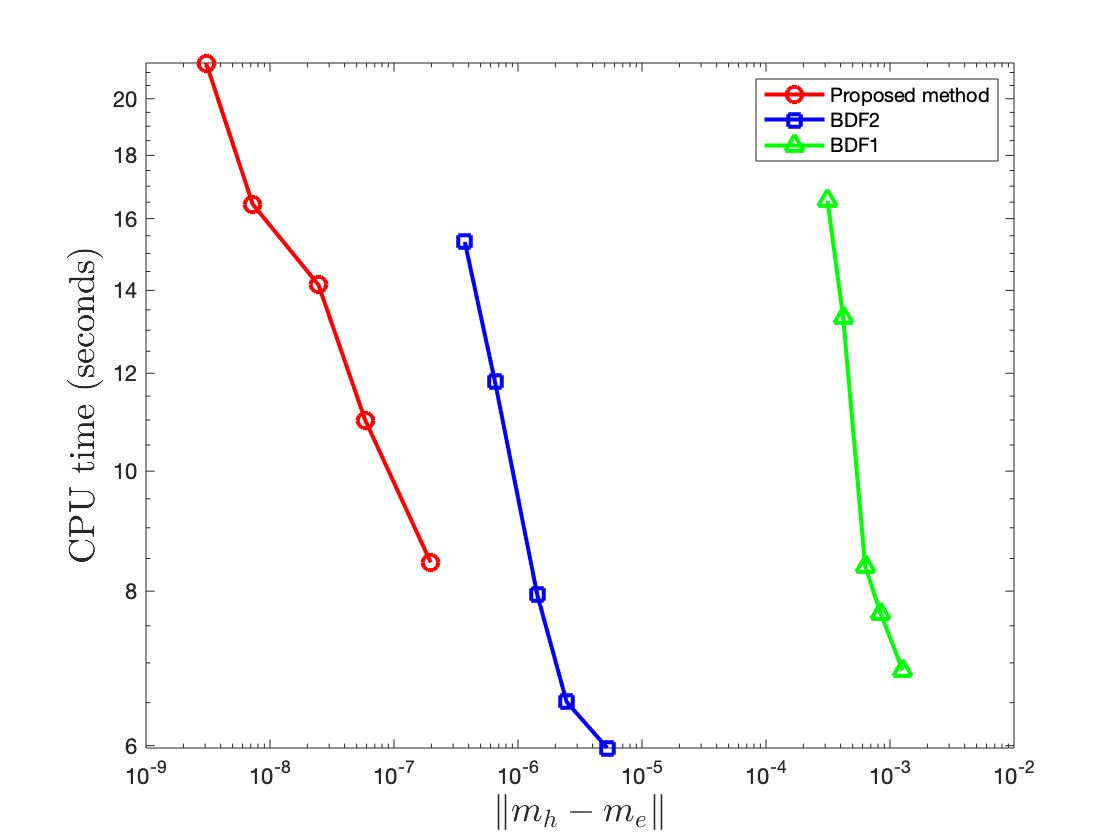}}
	\subfloat[Varying $k$ in 3D up to $T=0.1$]{\label{cputime_3D}\includegraphics[width=2.5in]{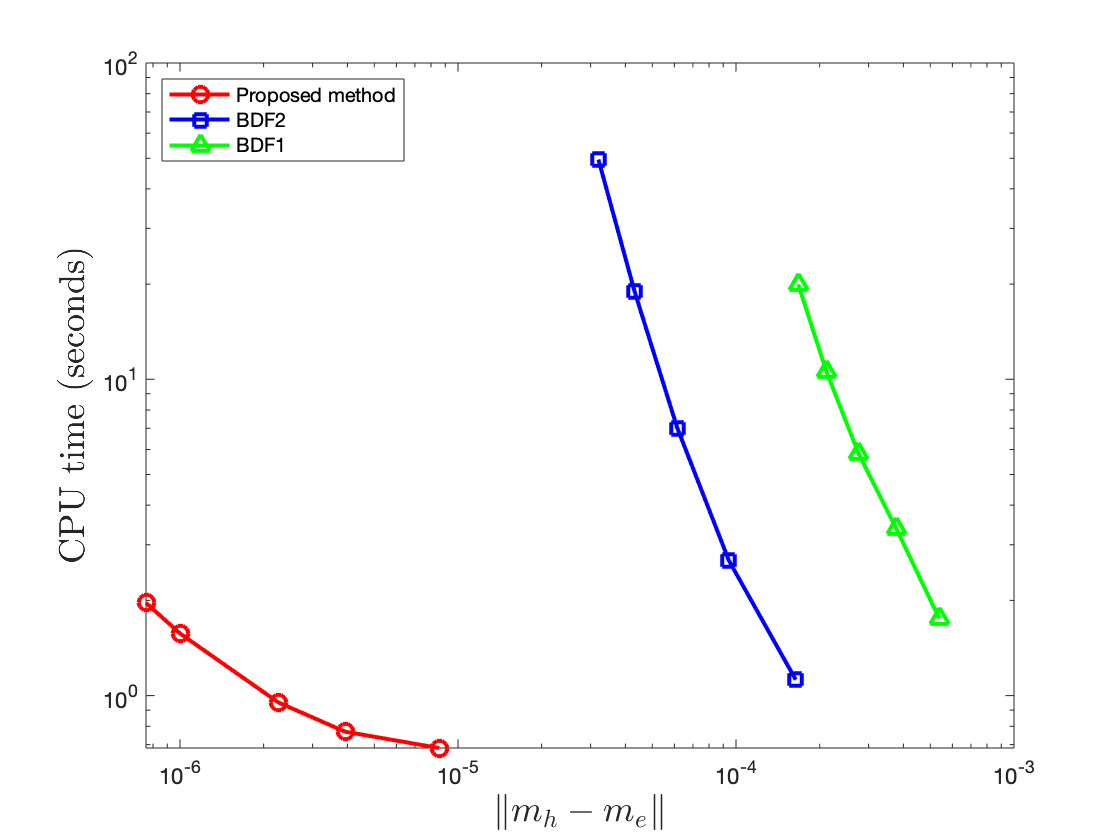}}
	\hspace{0.1in}
	\subfloat[Varying $h$ in 1D up to $T=1$]{\label{cputime_1D_space}\includegraphics[width=2.5in]{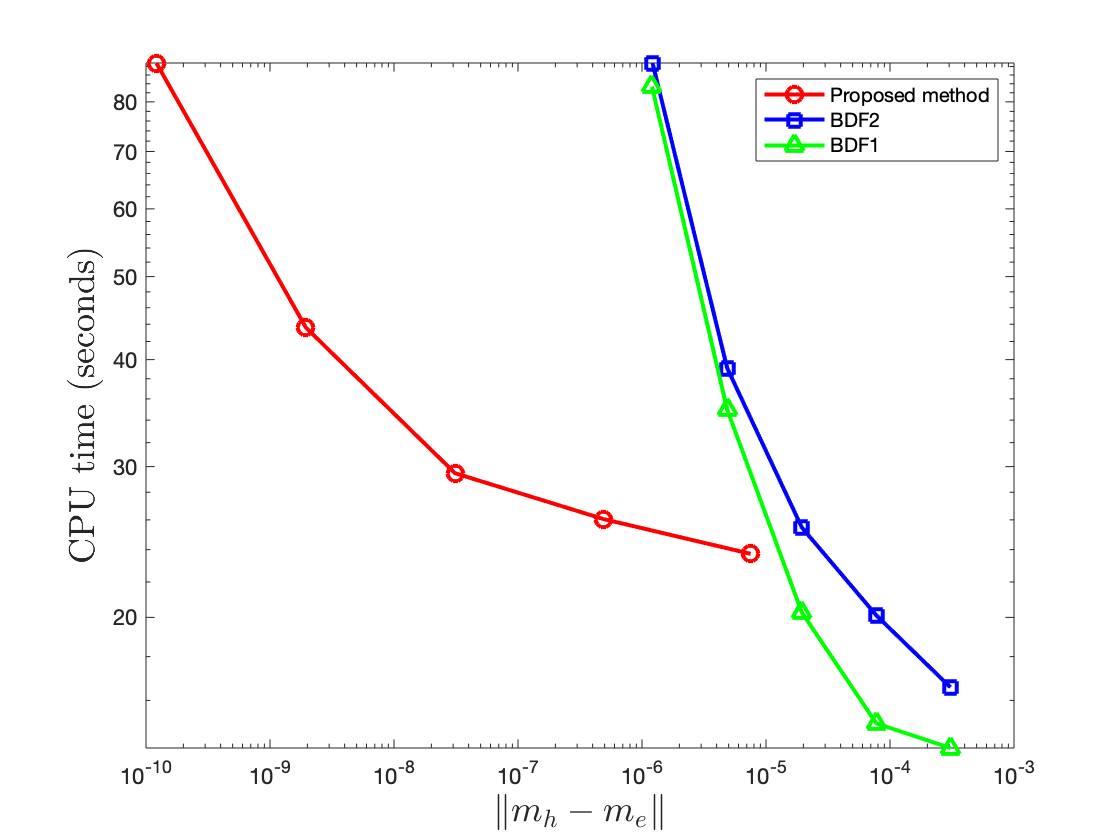}}
	\subfloat[Varying $h$ in 3D up to $T=1$]{\label{cputime_3D_space}\includegraphics[width=2.5in]{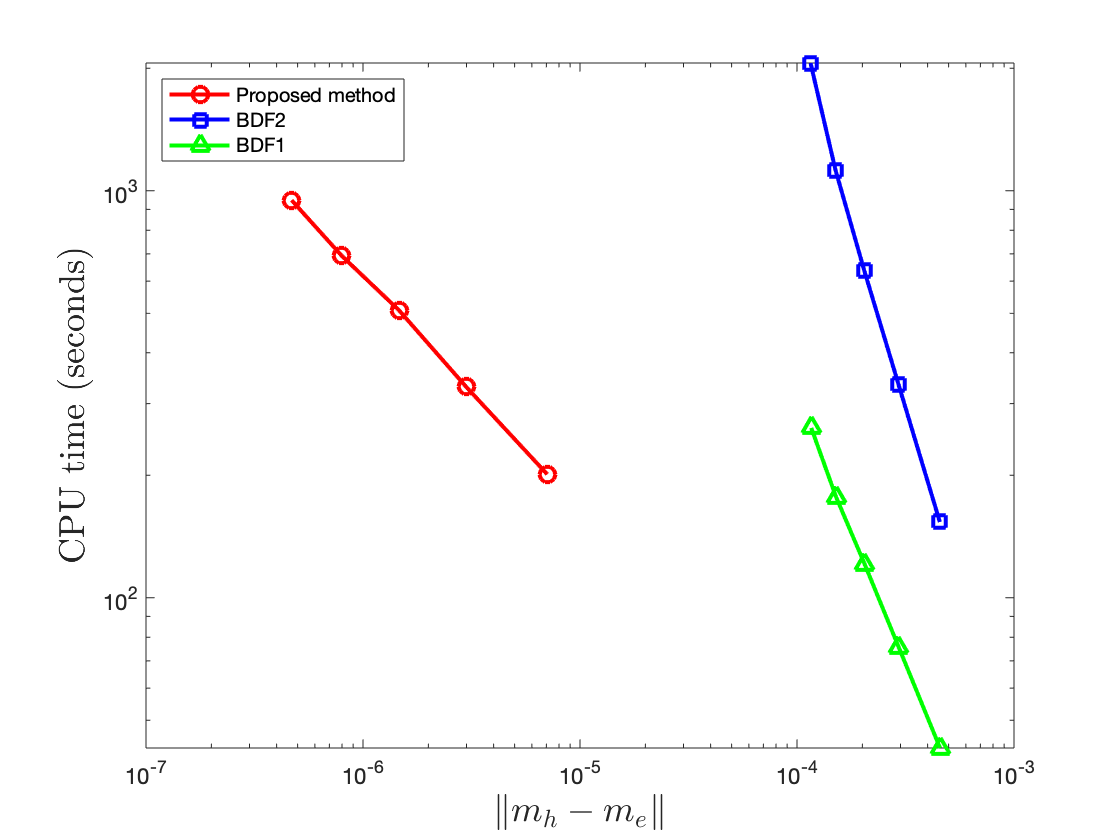}}
	\caption{CPU time needed to achieve the desired numerical accuracy, for the proposed method, the BDF2 and the BDF1, in both the 1D and 3D computations. The CPU time is recorded as a function of the approximation error by varying $k$ or $h$ independently. CPU time with varying $k$: proposed method $<$ BDF2 $<$ BDF1; CPU time with varying $h$: proposed method $<$ BDF1 $\lessapprox$ BDF2.}\label{cputime}
\end{figure}

To make a comparison in terms of the numerical efficiency, we plot the CPU time (in seconds) vs. the error norm $\|\m_h-\m_e\|_{\infty}$. In details, the CPU time is recorded as a function of the approximation error in \cref{cputime_1D} in 1D and in \cref{cputime_3D} in 3D, with a variation of $k$ and a fixed value of $h$. Similar plots are also displayed in \cref{cputime_1D_space} in 1D and \cref{cputime_3D_space} in 3D, with a variation of $h$ and a fixed value of $k$. In the case of a fixed spatial resolution $h$, the proposed method is significantly more efficient than the BDF1 and the BDF2 methods in both the 1D and 3D computations. The BDF2 is slightly more efficient than the BDF1, while such an advantage may vary for different values of $k$ and $h$. In the case of a fixed time step size $k$, the proposed method is more efficient than the BDF2 and BDF1, in both the 1D and 3D computations, and the cost of BDF1 is comparable to the
BDF2.

\subsection{Stability test with large damping parameters}
\begin{figure}[htbp]
	\centering
	\subfloat{\label{BDF3_alpha_1_ang}\includegraphics[width=1.3in]{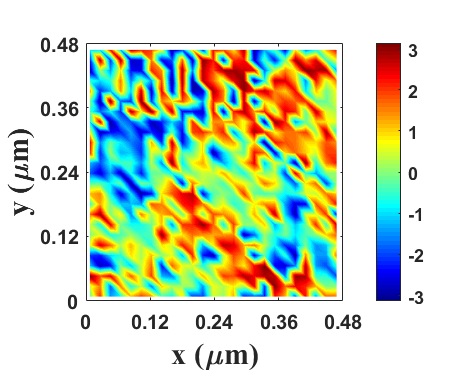}}
	\subfloat{\label{BDF3_alpha_5_ang}\includegraphics[width=1.3in]{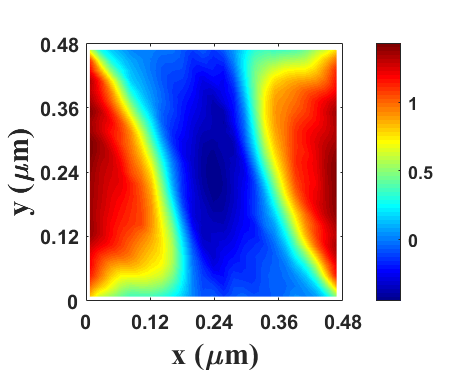}}
	\subfloat{\label{BDF3_alpha_10_ang}\includegraphics[width=1.3in]{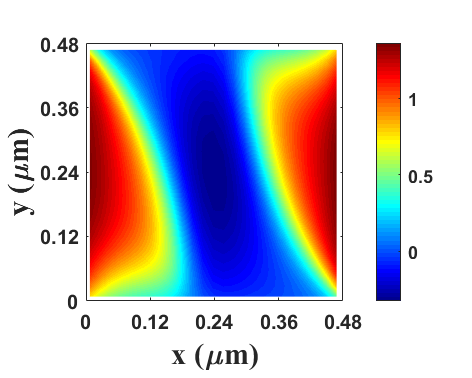}}
	\subfloat{\label{BDF3_alpha_40_ang}\includegraphics[width=1.3in]{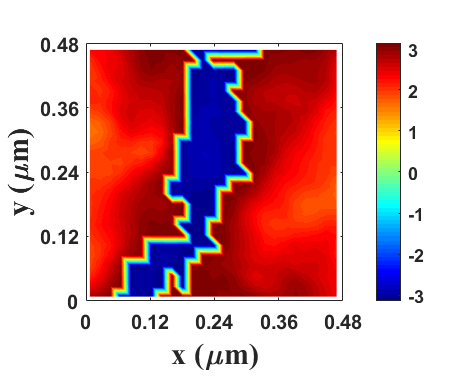}}
	\subfloat{\label{BDF3_alpha_100_ang}\includegraphics[width=1.3in]{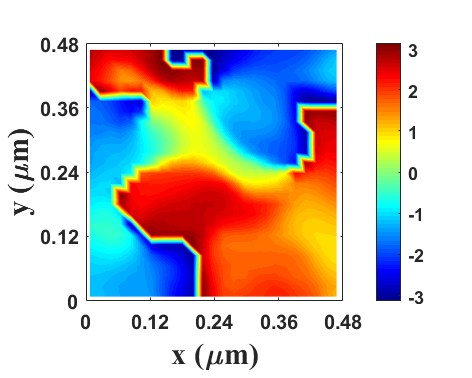}}
	\hspace{0.1in}
	\subfloat{\label{BDF2_alpha_1_ang}\includegraphics[width=1.3in]{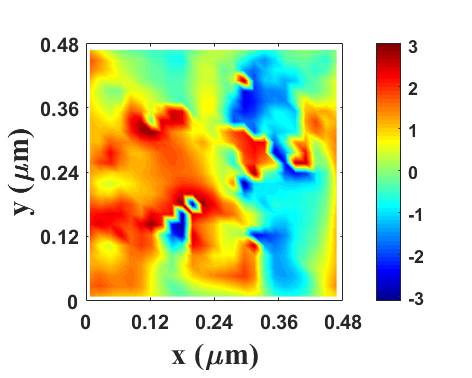}}
	\subfloat{\label{BDF2_alpha_5_ang}\includegraphics[width=1.3in]{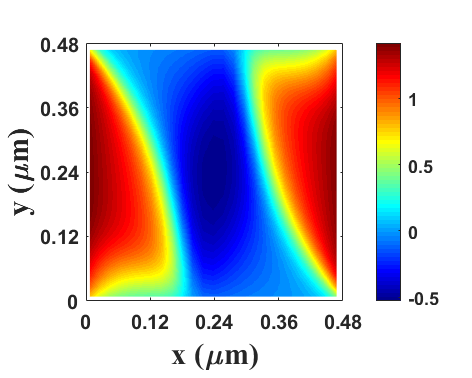}}
	\subfloat{\label{BDF2_alpha_10_ang}\includegraphics[width=1.3in]{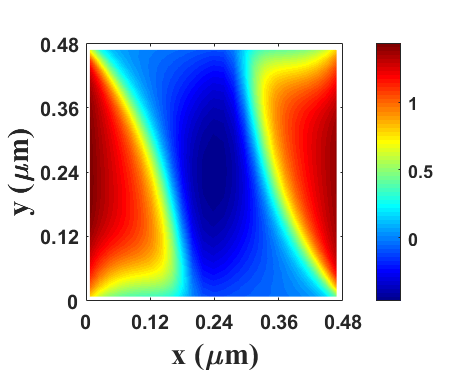}}
	\subfloat{\label{BDF2_alpha_40_ang}\includegraphics[width=1.3in]{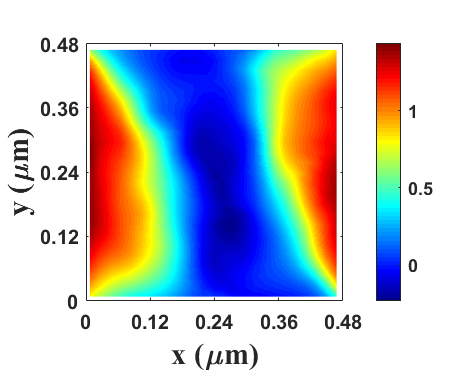}}
	\subfloat{\label{BDF2_alpha_100_ang}\includegraphics[width=1.3in]{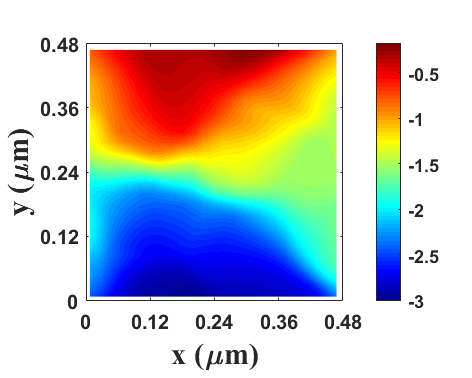}}
	\hspace{0.1in}
	\subfloat{\label{BDF1_alpha_1_ang}\includegraphics[width=1.3in]{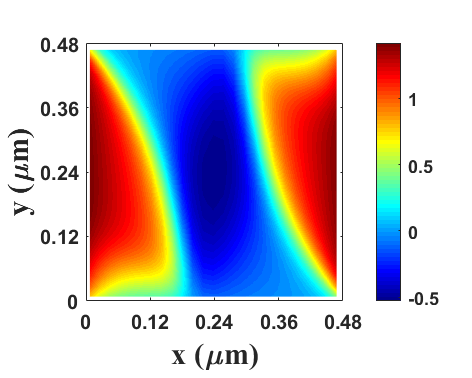}}
	\subfloat{\label{BDF1_alpha_5_ang}\includegraphics[width=1.3in]{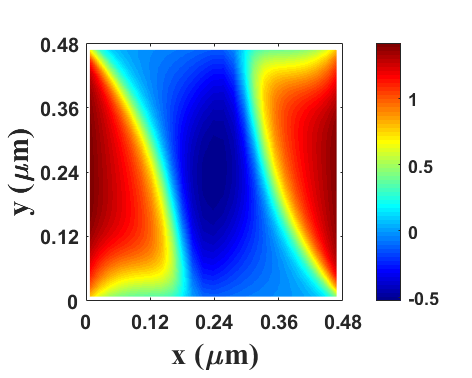}}
	\subfloat{\label{BDF1_alpha_10_ang}\includegraphics[width=1.3in]{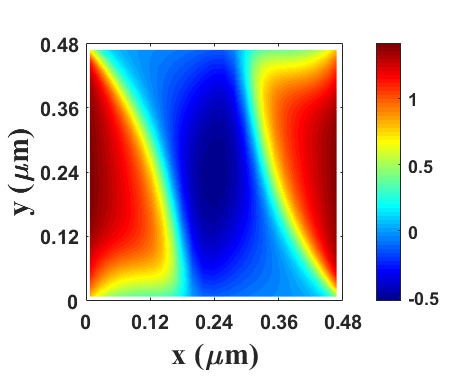}}
	\subfloat{\label{BDF1_alpha_40_ang}\includegraphics[width=1.3in]{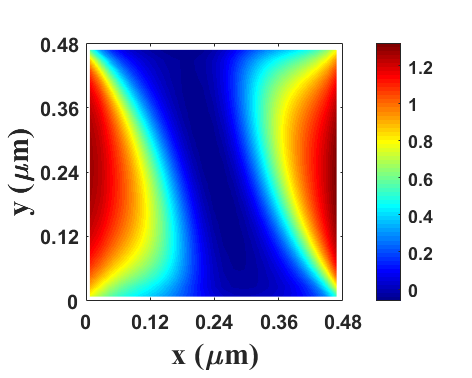}}
	\subfloat{\label{BDF1_alpha_100_ang}\includegraphics[width=1.3in]{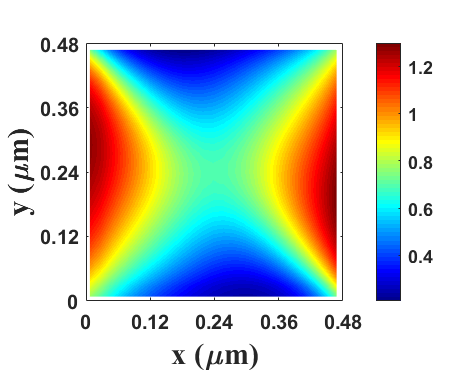}}
	\caption{Stable structures in the absence of magnetic field at $2\,$ns. The color denotes the angle between the first two components of the magnetization vector. Top: Proposed method; Middle: BDF2; Bottom: BDF1. From left to right: $\alpha=1,5,10,40,100$. $dt=1\;ps$. }\label{BDF3_BDF2_BDF1_alpha}
\end{figure}

\begin{figure}[htbp]
	\centering
	\subfloat{\label{BDF3_alpha_1_ang_v1}\includegraphics[width=1.3in]{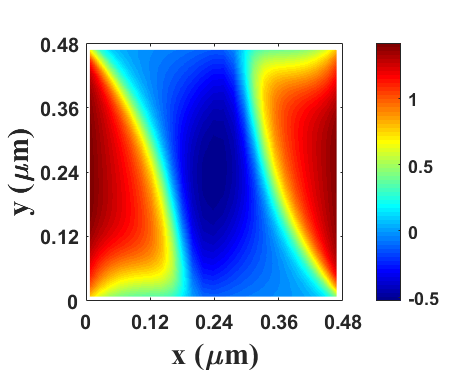}}
	\subfloat{\label{BDF3_alpha_5_ang_v1}\includegraphics[width=1.3in]{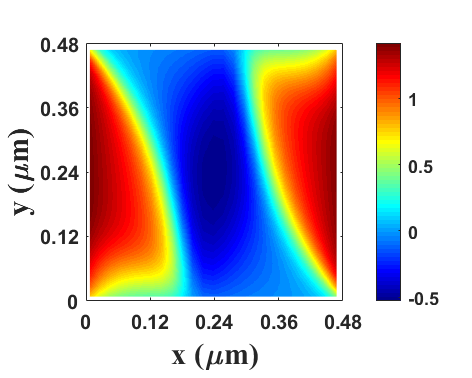}}
	\subfloat{\label{BDF3_alpha_10_ang_v1}\includegraphics[width=1.3in]{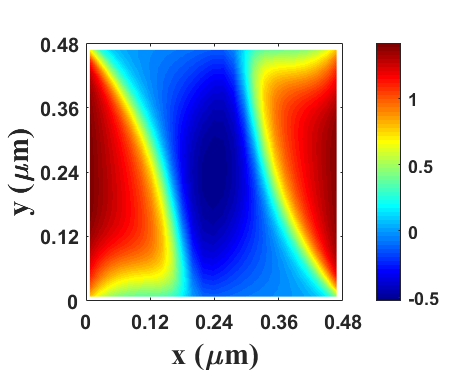}}
	\subfloat{\label{BDF3_alpha_40_ang_v1}\includegraphics[width=1.3in]{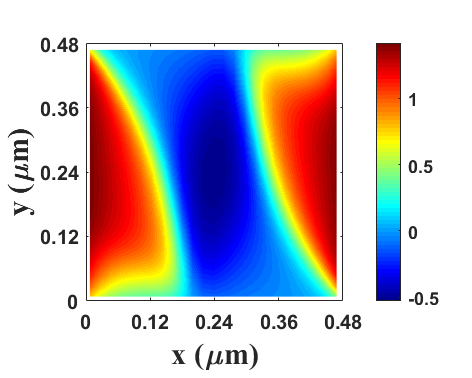}}
	\subfloat{\label{BDF3_alpha_100_ang_v1}\includegraphics[width=1.3in]{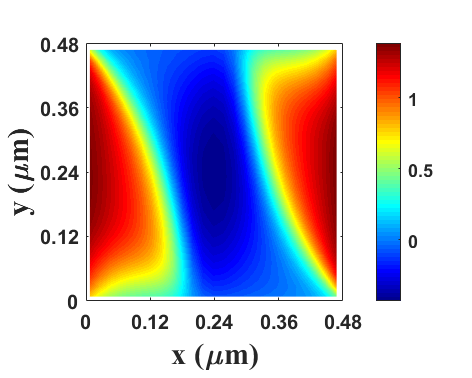}}
	\hspace{0.1in}
	\subfloat{\label{BDF2_alpha_1_ang_v1}\includegraphics[width=1.3in]{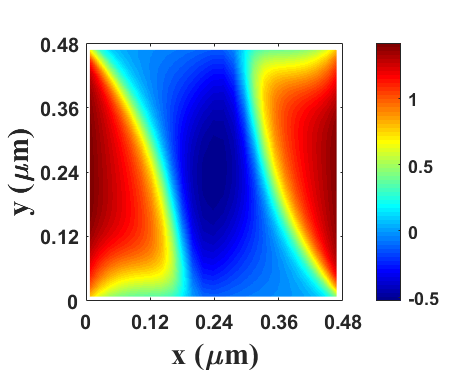}}
	\subfloat{\label{BDF2_alpha_5_ang_v1}\includegraphics[width=1.3in]{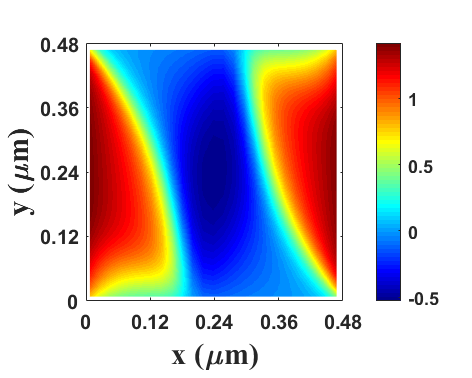}}
	\subfloat{\label{BDF2_alpha_10_ang_v1}\includegraphics[width=1.3in]{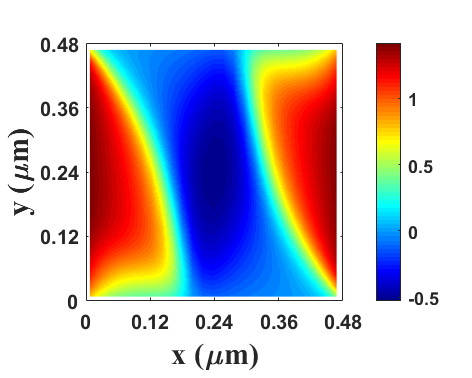}}
	\subfloat{\label{BDF2_alpha_40_ang_v1}\includegraphics[width=1.3in]{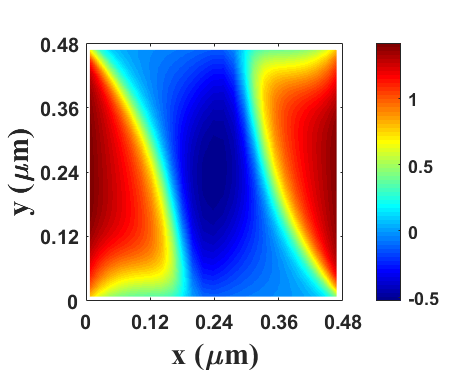}}
	\subfloat{\label{BDF2_alpha_100_ang_v1}\includegraphics[width=1.3in]{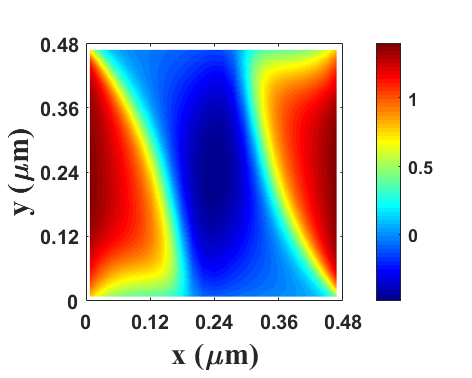}}
	\hspace{0.1in}
	\subfloat{\label{BDF1_alpha_1_ang_v1}\includegraphics[width=1.3in]{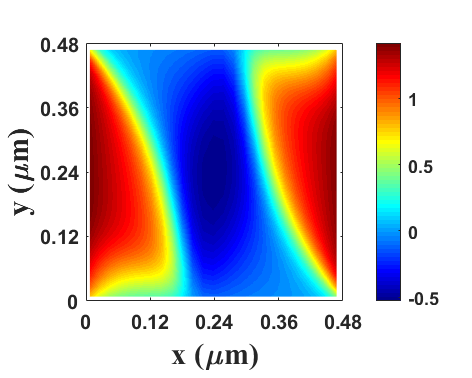}}
	\subfloat{\label{BDF1_alpha_5_ang_v1}\includegraphics[width=1.3in]{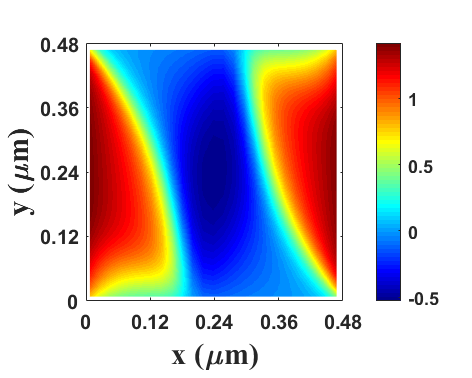}}
	\subfloat{\label{BDF1_alpha_10_ang_v1}\includegraphics[width=1.3in]{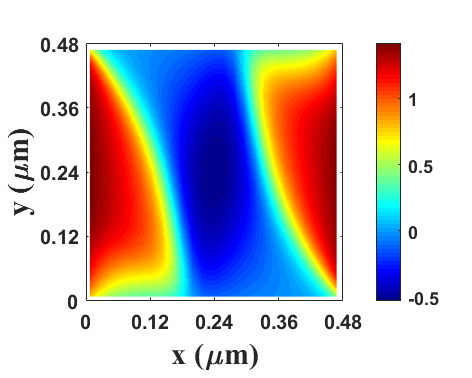}}
	\subfloat{\label{BDF1_alpha_40_ang_v1}\includegraphics[width=1.3in]{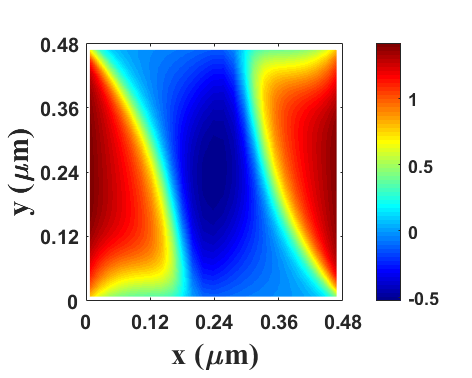}}
	\subfloat{\label{BDF1_alpha_100_ang_v1}\includegraphics[width=1.3in]{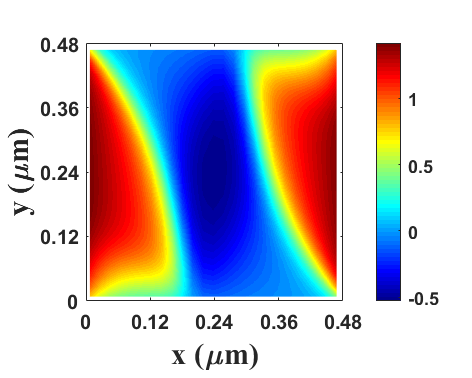}}
	\caption{Stable structures in the absence of magnetic field at $2\,$ns. The color denotes the angle between the first two components of the magnetization vector. Top: Proposed method; Middle: BDF2; Bottom: BDF1. From left to right: $\alpha=1,5,10,40,100$. $dt=0.1\;ps$. }\label{BDF3_BDF2_BDF1_alpha_v1}
\end{figure}

To check the numerical stability of these three methods in the practical simulations of micromagnetics with large damping parameters, we consider a thin film of size $480\times480\times20\,\textrm{nm}^3$ with grid points $100\times100\times4$. The temporal step-size is taken as $k=1\,$ps and $k=0.1\,$ps. A uniform state along the $x$ direction is set to be the initial magnetization and the external magnetic field is set to be $0$. Five different damping parameters, $\alpha=1,5,10,40,100$, are tested with stable magnetization profiles shown in \cref{BDF3_BDF2_BDF1_alpha} and \cref{BDF3_BDF2_BDF1_alpha_v1}. In particular, the following observations are made. 
\begin{itemize}
	\item All three methods are stable for small time step size $k=0.1\;$ps;
	\item All three methods are stable for moderately large $\alpha$ ($\alpha=5,10$) with slightly large $k=1\;$ps;
	\item As the order of the BDF method increases, its stability and applicability with respect to $\alpha$ (damping parameter) become narrower. For small values of $\alpha$, both the BDF2 and BDF3 methods are unstable; for larger values of $\alpha$, the BDF3 method is unstable.
\end{itemize}
In fact, a preliminary theoretical analysis reveals that, an optimal rate convergence estimate of the proposed method could be theoretically justified for $\alpha>7$ (however, BDF2 theoretically justified for $\alpha>3$). Meanwhile, extensive numerical experiments have implied that $\alpha>3$ (BDF2 with $\alpha>1$) is sufficient to ensure the numerical stability in the practical computations.

\begin{figure}[htbp]
	\centering
	\subfloat[Proposed]{\label{energy_BDF3}\includegraphics[width=1.8in]{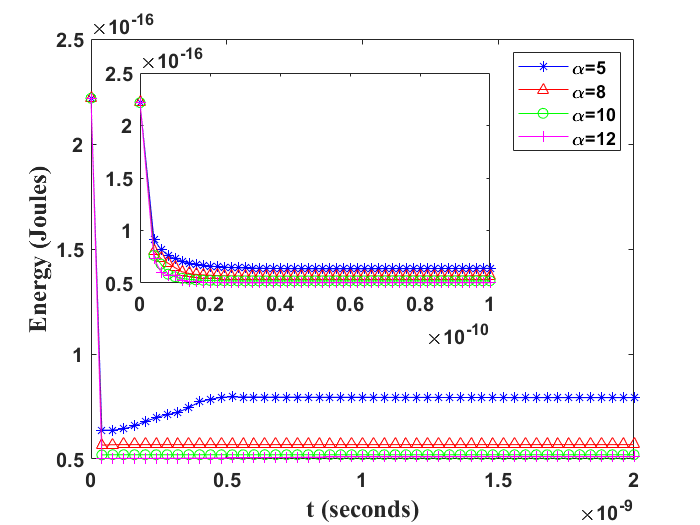}}
	\subfloat[BDF2]{\label{energy_BDF2}\includegraphics[width=1.8in]{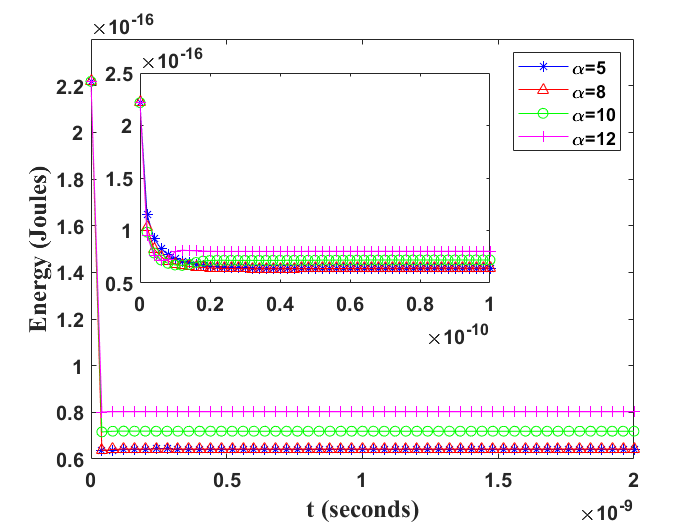}}
	\subfloat[BDF1]{\label{energy_BDF1}\includegraphics[width=1.8in]{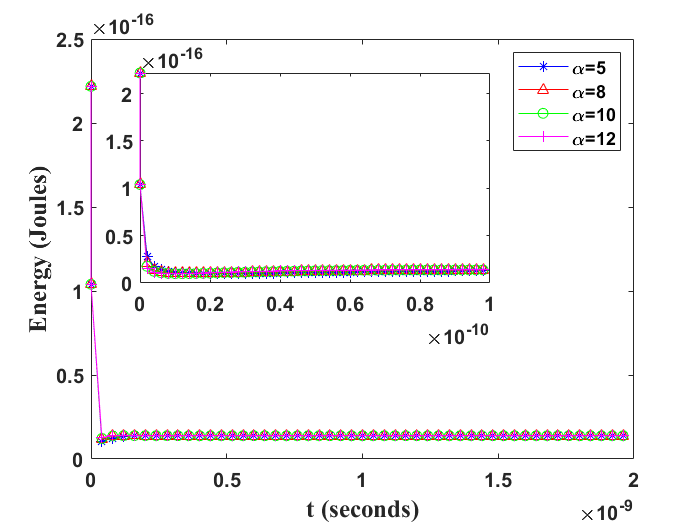}}
	\hspace{0.1in}
	\subfloat[Proposed]{\label{energy_BDF3_v1}\includegraphics[width=1.8in]{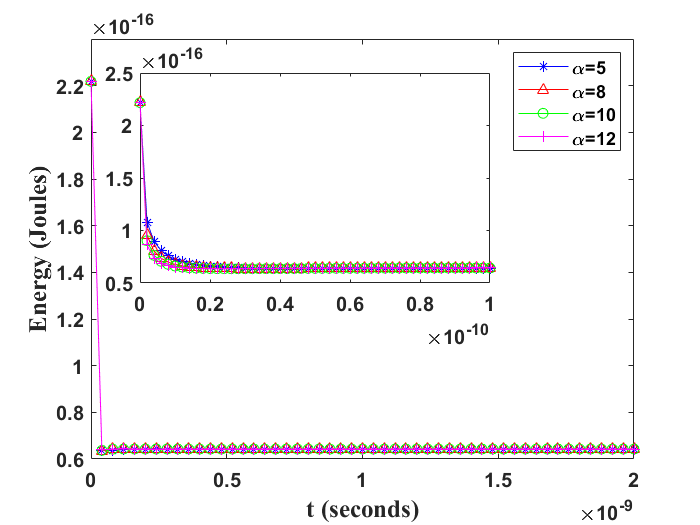}}
	\subfloat[BDF2]{\label{energy_BDF2_v1}\includegraphics[width=1.8in]{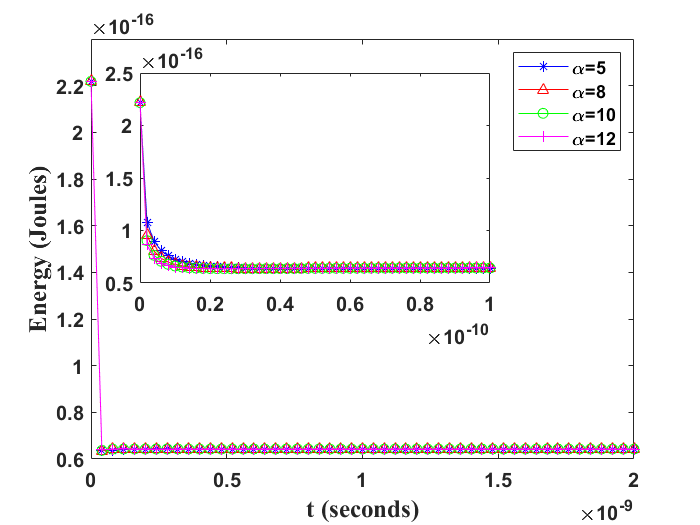}}
	\subfloat[BDF1]{\label{energy_BDF1_v1}\includegraphics[width=1.8in]{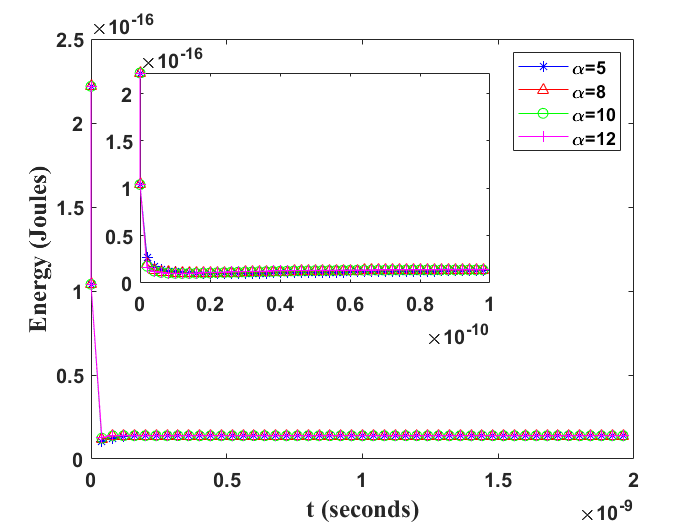}}
	\caption{Energy evolution curves of three numerical methods, with different damping constants, $\alpha=5,8,10,12$, up to $t=2\,$ns  in the absence of external magnetic field. Left: Proposed numerical method; Middle: BDF2; Right: BDF1. One common feature is that the energy dissipation rate is faster for larger $\alpha$, which is physically reasonable. Top row: $\Delta t=1\;ps$; Bottom row: $\Delta t=0.1\;ps$.}\label{energy_decay}
\end{figure}

Under the same setup outlined above, we investigate the energy dissipation of the proposed method, the BDF2, and the BDF1. The stable state is attainable at $t=2\,\textrm{ns}$ with $k=1\;$ps and $k=0.1\;$ps. The energy evolution curves of different numerical methods with different damping parameters, $\alpha=5,8,10,12$, are displayed in \cref{energy_decay}. One common feature is that the energy dissipation rate turns out to be faster for larger $\alpha$, in all three schemes. Meanwhile, a theoretical derivation also reveals that the energy dissipation rate in the LLG equation \eqref{eq-5} depends on $\alpha$, and a larger $\alpha$ leads to a faster energy dissipation rate. Therefore, the numerical results generated by all these three numerical methods have made a nice agreement with the theoretical derivation. The energy dissipation is much more stable among three methods with small time step size $k=0.1\;$ps. The convergent energy at $t=2\;$ns for BDF2 increases as $\alpha$ larger, however that of BDF1 and BDF3  decreases as $\alpha$ larger.

\begin{figure}[htbp]
	\centering
	\subfloat[$\alpha=5$]{\label{alpha_5}\includegraphics[width=2.5in]{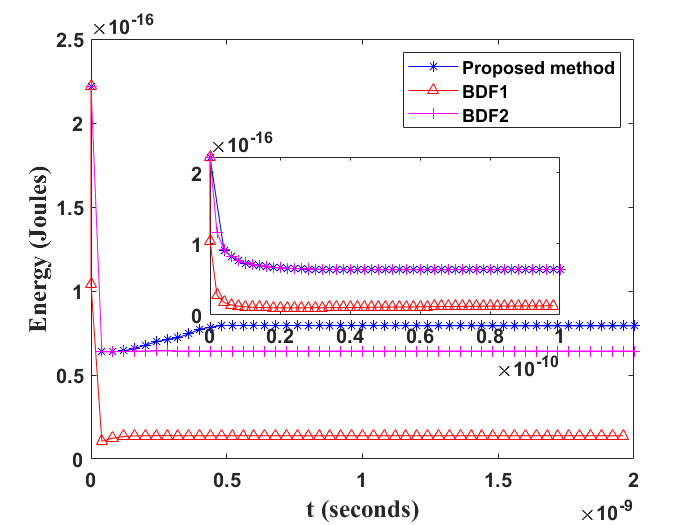}}
	\subfloat[$\alpha=8$]{\label{alpha_8}\includegraphics[width=2.5in]{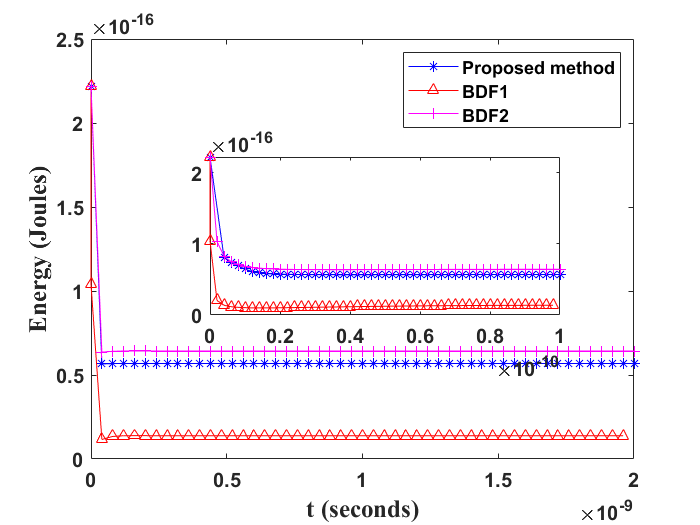}}
	\hspace{0.1in} 
	\subfloat[$\alpha=10$]{\label{alpha_10}\includegraphics[width=2.5in]{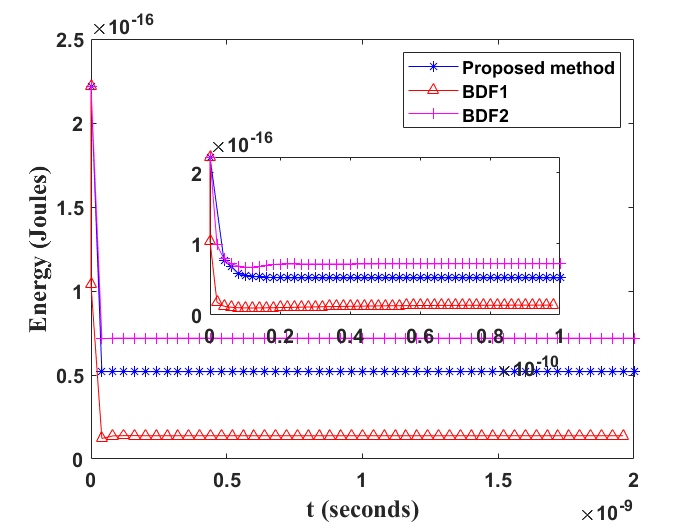}}
	\subfloat[$\alpha=12$]{\label{alpha_12}\includegraphics[width=2.5in]{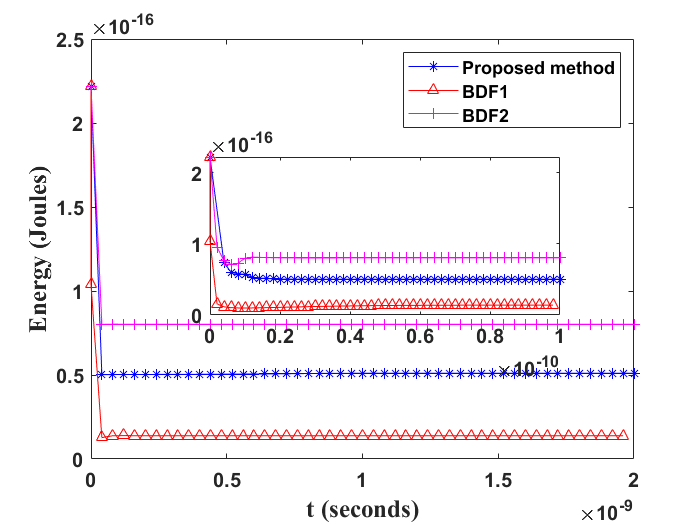}}
	\caption{Energy evolution curves in terms of time, for the numerical results created by three numerical methods up to $t=2\,$ns with $k=1\;ps$ in the absence of external magnetic field for (a) $\alpha=5$, (b) $\alpha=8$, (c) $\alpha=10$, and (d) $\alpha=12$. In a relatively short period of time, as the order of the BDF method increases, the energy decays more slowly. As $\alpha$ increases, the BDF1 method reaches the lowest energy at equilibrium, followed by the BDF3 method, and the BDF2 method has the highest energy. At smaller $\alpha$ values, the energy curves of BDF3 and BDF2 are basically the same. }\label{energy_decay_alpha}
\end{figure}

\begin{figure}[htbp]
	\centering
	\subfloat[$\alpha=5$]{\label{alpha_5_v1}\includegraphics[width=2.5in]{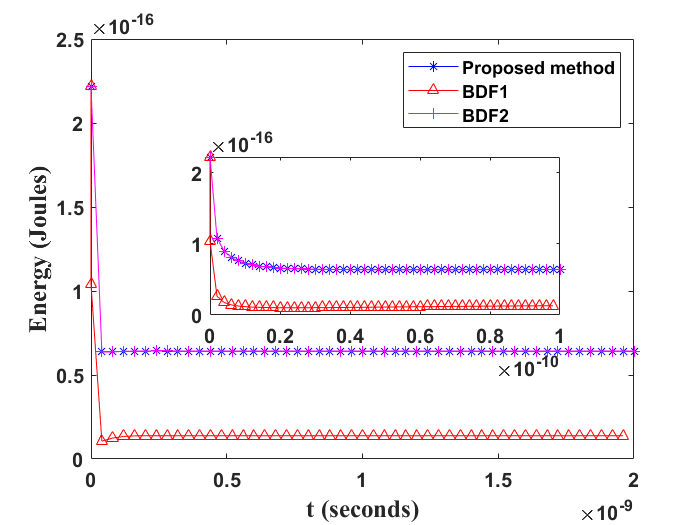}}
	\subfloat[$\alpha=8$]{\label{alpha_8_v1}\includegraphics[width=2.5in]{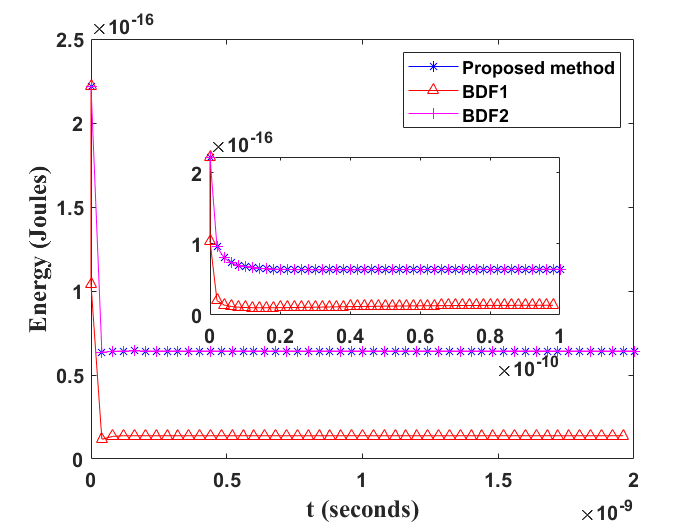}}
	\hspace{0.1in} 
	\subfloat[$\alpha=10$]{\label{alpha_10_v1}\includegraphics[width=2.5in]{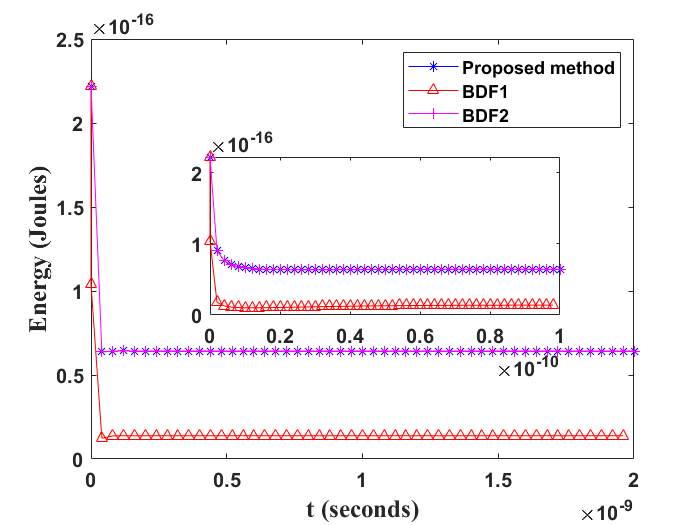}}
	\subfloat[$\alpha=12$]{\label{alpha_12_v1}\includegraphics[width=2.5in]{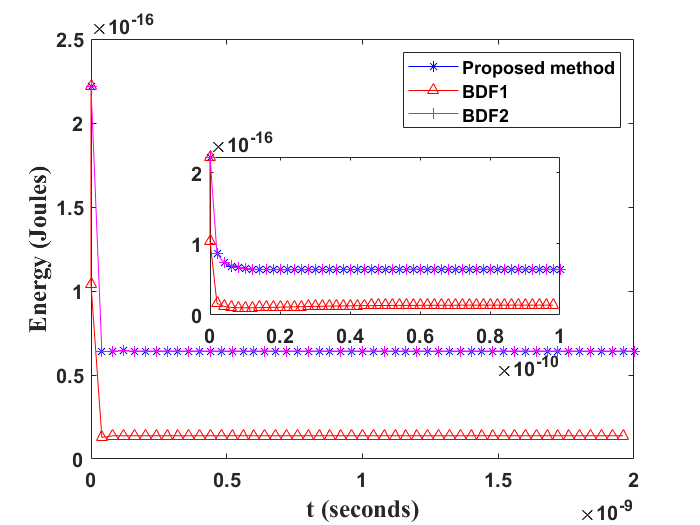}}
	\caption{Energy evolution curves in terms of time, for the numerical results created by three numerical methods up to $t=2\,$ns with $k=0.1\;ps$ in the absence of external magnetic field for (a) $\alpha=5$, (b) $\alpha=8$, (c) $\alpha=10$, and (d) $\alpha=12$. The energy dissipation pattern of the proposed method is consistent with the BDF2 method, and the BDF1 has a slightly different energy dissipation pattern from the other two methods and attains a lower level energy dissipation. The energy decay appears to be more stable.}\label{energy_decay_alpha_v1}
\end{figure}

Meanwhile, we choose the same sequence of values for $\alpha$, and display the energy evolution curves in terms of time up to $T=2\,$ns in \cref{energy_decay_alpha} and \cref{energy_decay_alpha_v1}. It is found that with $k=1\;$ps, in a relatively short period of time, as the order of the BDF method increases, the energy decays more slowly. As $\alpha$ increases, the BDF1 method reaches the lowest energy at equilibrium, followed by the BDF3 method, and the BDF2 method has the highest energy. At smaller $\alpha$ values, the energy curves of BDF3 and BDF2 are basically the same. For small time step size $k=0.1\;$ps, the energy dissipation pattern of the proposed method is consistent with the BDF2 method, and the BDF1 has a slightly different energy dissipation pattern from the other two methods and attains a lower level energy dissipation. The energy decay appears to be more stable.

\subsection{Domain wall motion}
\begin{figure}[htbp]
	\centering
	\subfloat[Magnetization for initial state]{\label{NeelWall_initial_mag}\includegraphics[width=2.8in]{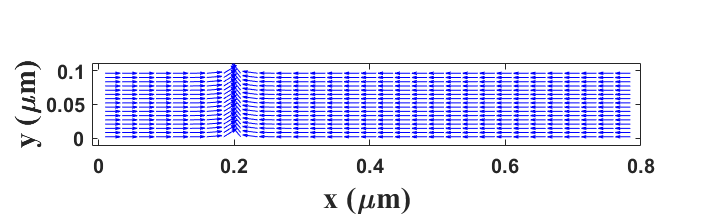}}
	\subfloat[Magnetization for initial state]{\label{NeelWall_initial_mag_v1}\includegraphics[width=2.8in]{BDF3_B_alpha_5_He_5mT_initial_dt_1ps.png}}
		\hspace{0.1in}
	\subfloat[Magnetization with $\alpha=5$  at $1.6\,$ns]{\label{NeelWall_alpha_5_2_mag}\includegraphics[width=2.8in]{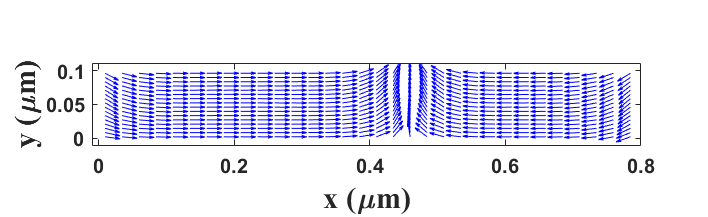}}
	\subfloat[Magnetization with $\alpha=5$  at $1.6\,$ns]{\label{NeelWall_alpha_5_2_mag_v1}\includegraphics[width=2.8in]{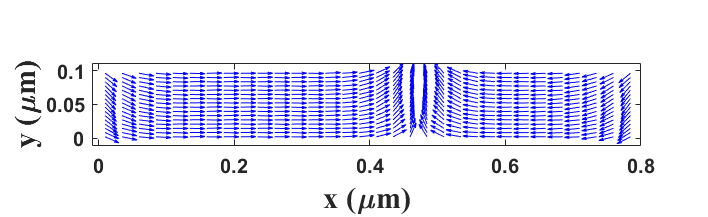}}
		\hspace{0.1in}
	\subfloat[Magnetization with $\alpha=8$ at $1.6\,$ns]{\label{NeelWall_alpha_8_2_mag}\includegraphics[width=2.8in]{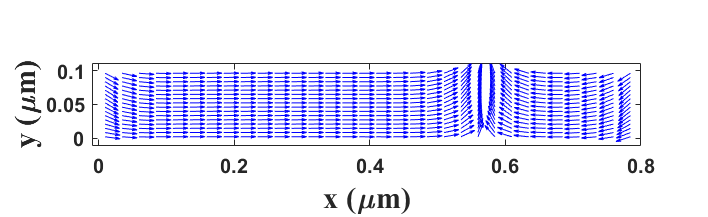}}
	\subfloat[Magnetization with $\alpha=8$ at $1.6\,$ns]{\label{NeelWall_alpha_8_2_mag_v1}\includegraphics[width=2.8in]{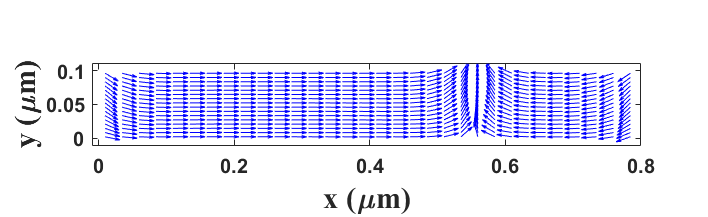}}
		\hspace{0.1in}
	\subfloat[Magnetization with $\alpha=10$ at $1.6\,$ns]{\label{NeelWall_alpha_10_2_mag}\includegraphics[width=2.8in]{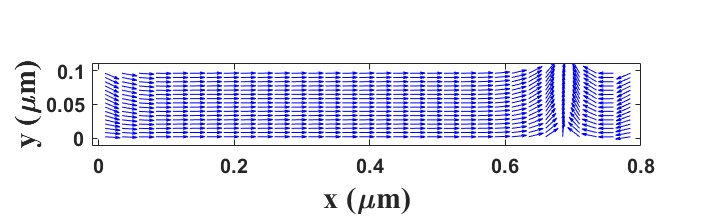}}
		\subfloat[Magnetization with $\alpha=10$ at $1.6\,$ns]{\label{NeelWall_alpha_10_2_mag_v1}\includegraphics[width=2.8in]{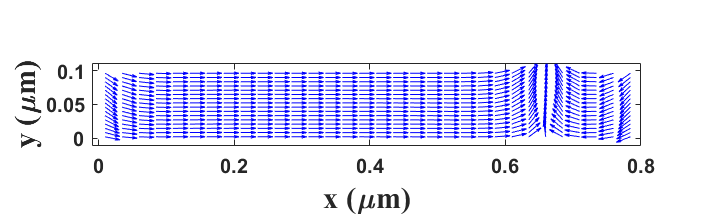}}
	\caption{Magnetization profiles of Ne\'{e}l wall motion in the presence of a magnetic field $\h_e=5\,$mT, with $\alpha = 5,8,10$ at $1.6\,$ns with $\Delta t=1\;$ps for the left panel and $\Delta t=0.1\;$ps for the right panel for the proposed numerical method. The in-plane arrow denotes the first two components of the magnetization vector. The wall moves faster for larger values of $\alpha$ and its velocity depends linearly on $\alpha$.}\label{NeelWall_alpha_2ns}
\end{figure}

\begin{figure}[htbp]
	\centering
	\subfloat{\label{velocity_alpha_varied_He}\includegraphics[width=2.8in]{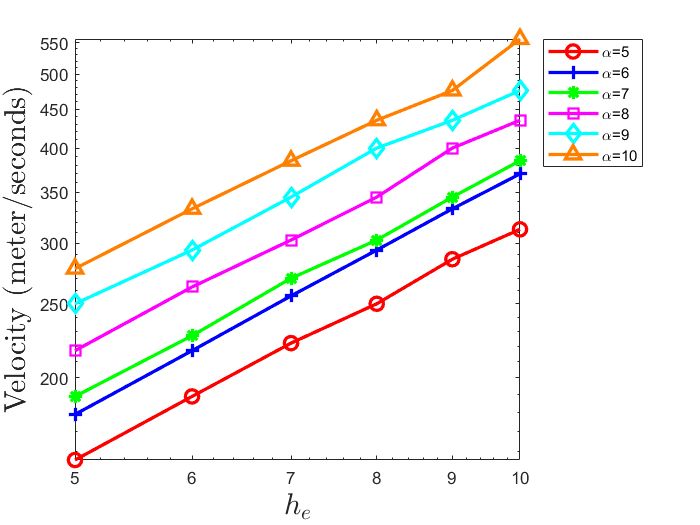}}
	\subfloat{\label{velocity_He_fixed_alpha}\includegraphics[width=2.8in]{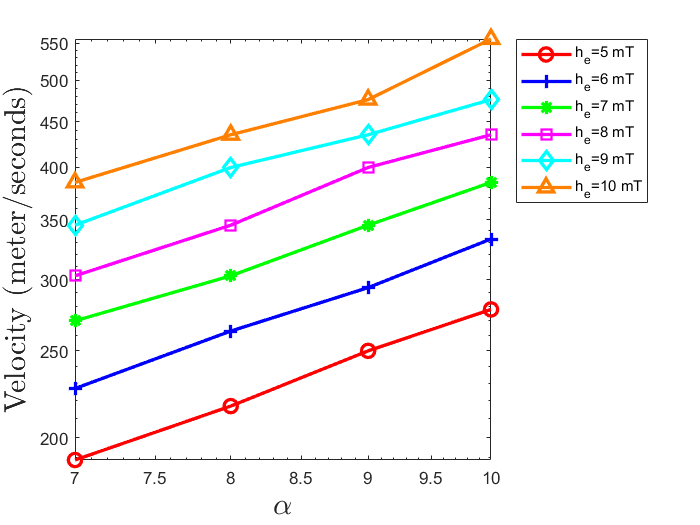}}
	\caption{Linear dependence of the wall velocity with respect to the damping parameter $\alpha$ (left) and the external magnetic field $\h_e$ (right).}\label{velocity_alpha_He}
\end{figure}

\begin{figure}[htbp]
	\centering
	\subfloat{\label{velocity_alpha_varied_He_v1}\includegraphics[width=2.8in]{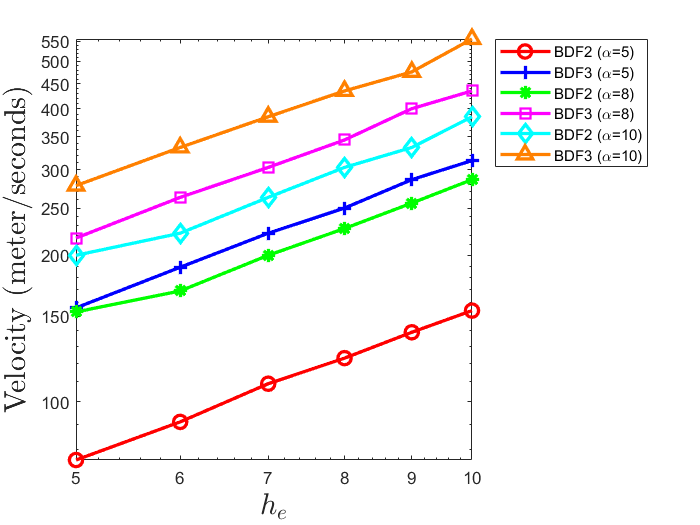}}
	\subfloat{\label{velocity_He_fixed_alpha_v1}\includegraphics[width=2.8in]{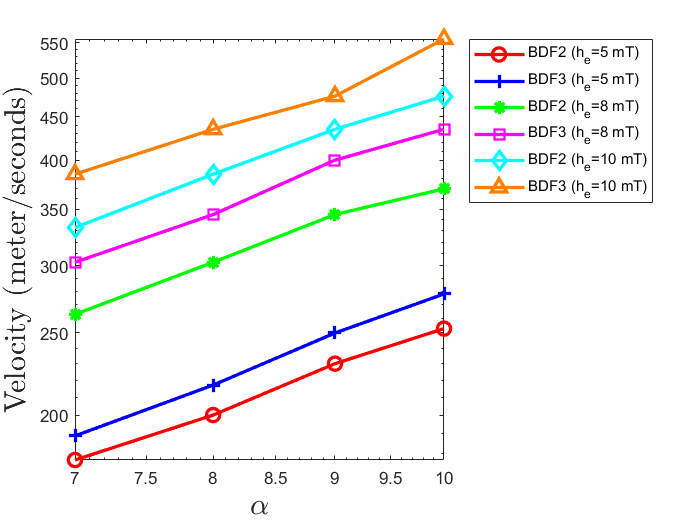}}
	\caption{Comparison of the wall velocity for BDF2 and BDF3 with respect to the damping parameter $\alpha$ (left) and the external magnetic field $\h_e$ (right).}\label{velocity_alpha_He_BDF2_BDF3}
\end{figure}

A Ne\'el wall is initialized in a nanostrip of size $800\times100\times4\,\textrm{nm}^3$ with grid points $128\times64\times4$. An external magnetic field of $\h_e=5\,$mT is then applied along the positive $x$ direction and the domain wall dynamics is simulated up to $1.6\,$ns with $\alpha=5,8,10$. The corresponding magnetization profiles are visualized in \cref{NeelWall_alpha_2ns} with $\Delta t=1\;$ps and $\Delta t=0.1\;$ps. Qualitatively, the domain wall moves faster as  the value of $\alpha$ increases. Quantitatively, the corresponding dependence is found to be linear; see \cref{velocity_alpha_He}. This conclusion is consistent with the previous conclusion obtained by simulating domain wall dynamics using the BDF2 method in \cite{CaiChenWangXie2022}. The slopes fitted by the least-squares method in terms of $\alpha$ and $\h_e$ are recorded in \cref{tab-3}. We can observe from \cref{velocity_alpha_He_BDF2_BDF3} that the velocity of domain wall motion obtained by the BDF3 method is faster than that obtained by the BDF2 method.

\begin{table}[htbp]
	\centering
	{\caption{Linear dependence of the domain wall velocity $V$ in terms of the external magnetic field $\h_e$ and the damping parameter $\alpha$.} \label{tab-3} }{
		\begin{tabular}{c|c|c|c|c|c|c|c}
			\hline 
			\diagbox{$\alpha$}{$V$ (m/s)}{$\h_e(\textrm{mT})$}&5&6&7 &8 & 9&10 & Slope\\
			\hline
			5&  156& 189& 222 &  250 & 286& 313 &1.007\\
			6&  179 & 217 & 256 & 294 & 333 & 370 &1.050\\
			7&  189&  227&  270 & 303 & 345& 385& 1.024\\
			8& 217 &  263 & 303 & 345 & 400 &435 &1.010 \\
			9& 250 & 294 & 345 &  400 &435 &476 &0.946 \\
			10& 278 & 333& 385 & 435 &476 & 556&0.965 \\\hline
			Slope & 1.094 &1.061  &1.003  & 1.039 &0.887 &0.998 &--\\
			\hline 
		\end{tabular}
	}
\end{table}

\section{Conclusions}
\label{sec:conclusions}


In the present study, a third-order accurate numerical method is proposed for solving the Landau-Lifshitz-Gilbert (LLG) equation with large damping coefficients. For the sake of numerical convenience, the LLG system is reformulated such that its damping term is re-expressed as a harmonic mapping flow. This numerical scheme is constructed based on the third-order backward-differentiation formula (BDF3) for approximating the temporal derivative, integrated with an implicit treatment of the constant-coefficient diffusion term, and a fully explicit extrapolation approximation for the nonlinear terms—including the gyromagnetic term and the nonlinear component of the harmonic mapping flow. Owing to the presence of large damping coefficients, the proposed method is verified to exhibit unconditional stability.
Compared with the second-order backward-differentiation formula (BDF2) method (with second-order temporal accuracy) and the first-order backward-differentiation formula (BDF1) method (with first-order temporal accuracy), the proposed method achieves higher accuracy. To validate the accuracy and computational efficiency of the proposed numerical method, numerical results in one-dimensional (1D) and three-dimensional (3D) domains are presented. Furthermore, micromagnetic simulations implemented via the proposed method yield physically consistent structures and successfully capture the linear dependence of domain wall velocity on the damping coefficient. Consequently, the proposed method can be efficiently applied to challenging practical simulations of micromagnetics involving large damping coefficients. Notably, the domain wall motion velocity obtained using the proposed method is higher than that derived from lower-order methods.

\section*{Acknowledgments}
This work is supported in part by the grants NSF DMS-2012669 (C.~Wang), and Jiangsu Science
and Technology Programme-Fundamental Research Plan Fund, Research and Development Fund
of XJTLU (RDF-24-01-015) (C. Xie).

\bibliographystyle{amsplain}
\bibliography{references}

\end{document}